\documentclass[12pt]{article}
\usepackage[dvips]{graphicx}

\begin{document}
\begin{center}
{\large \bf Modeling the influence of TH1 and TH2 type cells in autoimmune diseases}
\end{center}

Key words: EAE,IDDM,Mathematical models,TH1,Network\\

\begin{center}
{\large \bf Yoram Louzoun{\scriptsize (1)},Henri Atlan{\scriptsize (2)},Irun. R. Cohen{\scriptsize (3)}}  
\end{center}

\textit{\scriptsize (1) Interdisciplinary Center for Neural Computation, Hebrew
University, Jerusalem Israel.  yoraml@alice.nc.huji.ac.il}{\scriptsize \par{}}{\scriptsize \par}

\textit{\scriptsize (2) Human Biology Research Center, Hadassah Hebrew University Hospital, Jerusalem Israel.}{\scriptsize \par{}}{\scriptsize \par}

\textit{\scriptsize (3) Department of Immunology, The Weizmann Institute of Science.}{\scriptsize \par{}}{\scriptsize \par}

\newpage

\begin{abstract}
A sharp TH1/TH2 dichotomy has often been used to define the effects of cytokines on autoimmune diseases.
However contradictory results in recent research indicate that the situation may be more complex. 
We build here a simple mathematical model aimed at settling the contradictions. 
The model is based on a neural network paradigm, and is applied using Partial Differential Equations (PDE). 
We show here that a TH1/TH2 paradigm is only an external view of a complex multivariate system. 
\end{abstract}

\newpage

\section{Introduction}

Immunologists generally agree that in the course of many autoimmune diseases, TH1 type cells have a damaging effect, whereas TH2 type cells have a beneficial effect \cite{R4}. 
Thus the disease is aggravated by a high ratio of TH1/TH2 T cells and cytokines. 
This view is contradicted by many experimental results showing an inverse effect \cite{R2}: TH2 type cells may enhance disease, if administrated at a certain time \cite{R5}, and TH1 type cells can stop the evolution of the disease, if given at the appropriate time \cite{R18,R6a}.
Other works show that gene knock-out mice with no INF $\gamma$ or TNF $\alpha$ still manage to produce the pathological features of a TH1-dependent autoimmune disease \cite{R7,R7b}.
Based on results from EAE (experimental autoimmune encephalomyelitis) \cite{R22a} and insulin dependent diabetes melitus (IDDM) \cite{R21,R21a} in Non Obese Diabetic (NOD) mice \cite{R9}, we propose a larger model explaining both the general agreement and the contradictory results. 
We suggest that the TH1 and TH2 cells are not the cause of the disease, but a marker for a more general steady state. 
To settle the contradictions, we use a numerical model and general mechanisms of safe information transmission \cite{R9a} to explain the discrepancy in the case of gene knock-out mice. 
We show the explanatory power of our model and its ability to produce the known features of TH1 type diseases.

\section{Model}

Our model is based on a neural network paradigm \cite{R9c}, where the overall situation is not described by the properties of any one cell type, but by the steady state of the network. 
We use PDE in order to get a continuous description of the network. Each cell in the network represents a class of cells with the same general properties, at least from the most relevant aspect of the TH1/TH2 dichotomy. 
The model is an ad hoc description of the dynamics, but most of the links we use have some experimental backing. We built the model with the smallest number of cell types needed to explain qualitatively the experimental results.

\subsection{Experimental results}
EAE is an autoimmune disease studied mainly in rats and in mice, and is used as a model for human multiple sclerosis \cite{R8,R8a} .
In EAE, myelin in the CNS becomes inflamed inducing paralysis and even death \cite{R10}. 
This disease is classically defined by a TH1/TH2 dichotomy in which TH1 T cells induce the disease whereas TH2 cells prevent it \cite{R4}.

Classical EAE experiments show that: 
\begin{itemize}
\item EAE can be induced either by adoptive transfer of CD4 T cells from a sick rodent, or by active immunization by an injection of a myelin antigen in adjuvant  to a healthy rodent \cite{R14}. We will show that these two procedures lead, in fact, to two separate pathways of disease induction.

\item Administration of free antigen arrests EAE \cite{R15}. We propose to show that antigen-induced inhibition could be due to the effect of regulatory anti-idiotypic (anti-id) suppressor cells.
 
\item The disease can be prevented by the continuous administration of IL-4, which is a TH2 type cytokine thought to push the immune system into a healthy TH2 dominated state  REF FROM EAE\cite{R16}. But IL-4 can also aggravate disease (see below).
\end{itemize}

However, recent research results appear to contradict the classical TH1/TH2 hypothesis:
\begin{itemize}
\item A small amount of IL-4 administrated during disease induction can block the suppressive effect of free antigen.\cite{R17}
This paradox will be explained by the effect of TH2 type cytokines on anti-id regulatory cells.

\item A short administration of a TH1 cytokine  can heal the disease \cite{R18,R6a}. This effect too can be attributed to an anti-id regulatory cell. 

\item  Gene knock-out mice, lacking one or some of the TH1 cytokines, can still develop EAE or IDDM \cite{R7,R7b}. This paradox will be explained by the nature of regulation.
\end{itemize}

These facts regarding EAE are combined in our analysis with results regarding the structure of id - anti-id cell interactions and the cytokine profile expressed during the evolution of IDDM in NOD mice \cite{R21,R21a} , and in the mouse model of EAE \cite{R22a,R8b}. 

\subsection{Cells and Cytokine Types}

Our model is composed of five cell types, three cytokine types and the antigen (Table 1). 
The first cell type is the naive CD4 T cell that can differentiate into the second and third cell types: TH1 and TH2 type CD4 T cells \cite{R22}.
The anti-id cells are grouped as a single cell type, the fourth cell type. A more detailed model should use naive anti-id cells and diverse anti-id cells secreting TH1 type cytokines and TH2 type cytokines. 
For the sake of simplicity, however, we group all anti-id cells into one cell type, and assume they have the same effect. Since we know from IDDM and EAE that anti-id cells secrete more TH1 cytokines than TH2 \cite{R21a,R22a} , we designate the anti-id population as TH1 type cells. The CD4 T id cells can present themselves directly to the anti-id cells \cite{R21a}.
The fifth type of cells are macrophages, which play additional roles as APC \cite{R23} and as cytokine-secreting cells. The antigen represents the quantity of antigen presented to the CD4 T cells. 
We assume that there is a large enough number of APCs so that the CD4 T cells can see the antigen. 

The three types of cytokines are:
\begin{itemize}
\item TH1 type cytokines that enhance the differentiation of naive cells into TH1 type cells. 
\item TH2 type cytokines that enhance the development of TH2 type cells. 
\item Cytokines secreted by macrophages that enhance the proliferation and differentiation of TH1 cells \cite{R24}. We shall ignore the effect of TH2 type cytokines secreted by macrophages (like IL-10) and the effects of B cells. We will group together all the  macrophage-derived  cytokines as type C cytokines.
\end{itemize}

\begin{table}
\begin{tabular}{|c|c|c|}\hline
Group&Cell Type&Notation\\
\hline
1&Naive cells&id - 0\\
\hline
2&Idiotypic TH1 cells&A\\
\hline
3&Idiotypic TH2 cells&B\\
\hline
4&Anti-id cells&anti-id\\
\hline
5&Macrophages&M$\Phi$\\
\hline
6&TH1 cytokines&Cyt A\\
\hline
7&TH2 cytokines&Cyt B\\
\hline
8&Cytokines secreted by macrophages&Cyt C\\
\hline
\end{tabular}
\caption{Cell types and cytokines. The cytokines are labeled according to the cell type producing them. We grouped cells as follows: TH1 id cells can include all types of lymphocytes producing TH1 cytokines. Anti-id cells are all cell types specifically reactive to the id cells. They can either be cytotoxic cells, or limit the proliferation rate of the id cells.}
\end{table}

\subsection{Model equations.}
We compute the time derivatives of each cell type concentration $x_i$ and of the cytokine concentration $y_i$, assuming both are affected by all other cells and cytokines. 

\begin{eqnarray}
x_i'=\Sigma w_{ij}x_i x_j  +\max(\Sigma w_{1ij}x_{origin i}x_j,0) -d_i x_i \\
y_i'=max(\Sigma w_{2ij} x_j,0)  -d_i y_i
\end{eqnarray}

$x_{origin i}$ is the origin of the cell number i (for example, the TH1 cells originate in the naive cells, the naive cells originate in themselves). 
The other cells and cytokines can either induce or inhibit the proliferation of this ($x_i$) cell type. The cytokines originate from the cells that secrete them, so that the right side of eq' 2 is not multiplied by the cytokine concentration. 
The secretion rate of the cytokines ($y_i$) by cells ($x_j$) is linear to the other cell and cytokine concentrations. For example, if TH2 type cytokines inhibit the production of TH1 cytokines, then we use a term in the equation for the TH1 cytokines $-w_{2ij}x(th_2)$, and we do not multiply the term by the TH1 cytokine concentration.

A cell's proliferation rate cannot be negative. Thus we define the proliferation rate of any cell to be greater than 0. But of course, the natural death rate ($-d_ix_i$) can lead to a decrease in the cell concentration. T
he cytokine production rate ($y_i'$) is assumed not to depend on the cytokine concentration, but on the number of cells secreting the cytokine.
The equations are solved numerically using a second order precision scheme \cite{R26}. 

We start every computation at a steady state 
\begin{equation}
x_i'=0  , y_i'=0
\end{equation}
Since the time constants in this process are short, we could assume that the system had already reached a steady state before the appearance of the autoimmune disease.

Macrophages do not reproduce during the response to the antigen but, rather, are attracted to the affected region \cite{R27}, so that, like the cytokines,their growth rate is not multiplied by their own concentration.

The differences between possible models are in the assigned weights $w_{ij}, w_{1ij}$ and $ w_{2ij}$ in eq' 1 and 2. But we have tried to use the most uniform definition of weights, and all the weights are of the same order of magnitude. If we had no special reason to believe that two weights of the same type were different, we fixed them at the same value (the effect of TH1 cells on TH1 cytokines and the effect of TH2 cells on TH2 cytokines, for example). This was done in order to limit the size of the parameter space. 

\subsection{Model description}

The structure of the model is simply a mathematical description of the findings in IDDM and EAE. The basic part of the model is the opposition between id cells and anti-id cells. The CD4 are the id cells and the naive CD4 population proliferates proportionally to the antigen concentration. The suppressors are the anti-id cells, and they proliferate proportionally to the concentration of id cells (Figure 1A).  

The second division of this model is into TH1 and TH2 id cells. Both the TH1 and TH2 cells originate from the naive id cell. The TH1 cell type is associated with a TH1 cytokine type and the TH2 cell type is associated with a TH2 cytokine type. We assume that each cell enhances its own cytokine type \cite{R27a} and inhibits the production of the opposing type. Each cytokine encourages the production of the appropriate cell type and inhibits the other type (Figure 1B). 

The macrophages have a double effect: They serve as APCs and activate naive CD4 cells \cite{R23}, and they secrete cytokines that enhance T cell differentiation into TH1 type cells \cite{R29} (and inhibit differentiation into TH2 cells). The macrophages are influenced by the TH1 cytokines, which activate them to migrate into the inflammatory site. Since macrophages respond to TH1 cytokine and are inhibited by TH2 cytokines, changes in the macrophage concentration depend on the differences between the concentrations of TH1 and TH2 cytokine types. 

The anti-id cells are of a TH1 type. This assumption is based on data from EAE and IDDM \cite{R21a,R22a}.
We thus postulate that the anti-id cells are inhibited by TH2 type cytokines. And finally, we assume that there is a basic concentration and differentiation rate of CD4 cells into TH1 and TH2; but the natural situation (when there are no cytokines involved)  tends toward the TH2 dominated state (a healthy situation). 
The above model is described in Figure 2

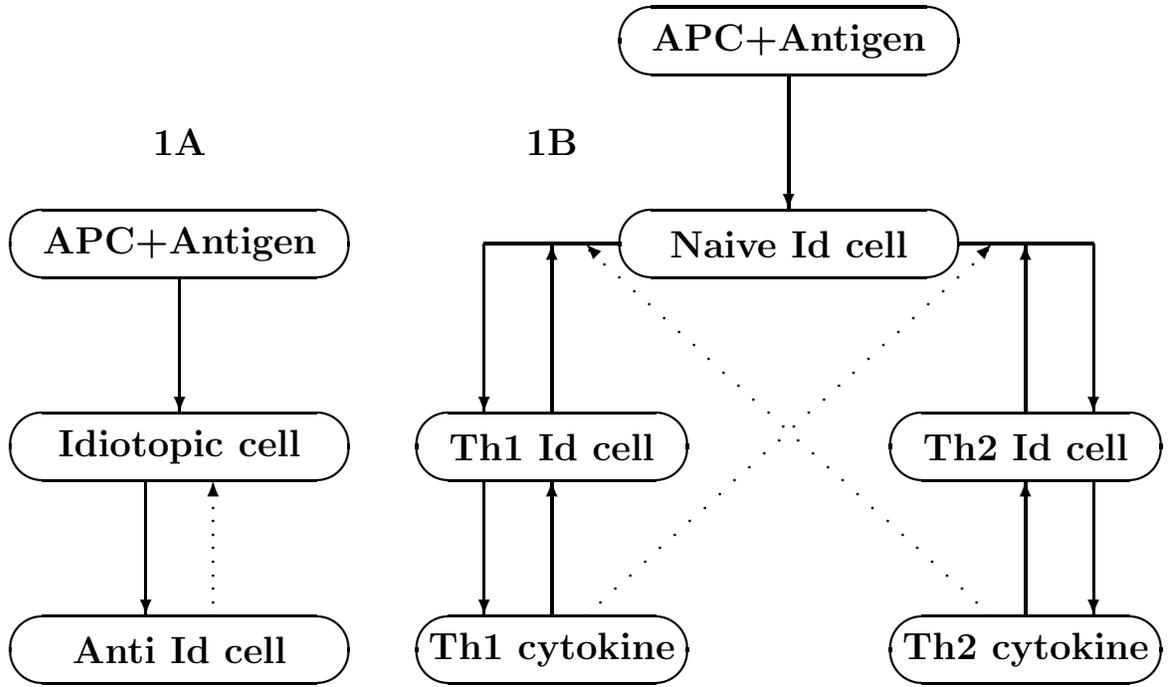
\begin{figure}
\center
\noindent
\setlength{\unitlength}{0.9cm}
\begin{picture}(20.0,11.0)
\thicklines
\put(2.0,8.5){\makebox(0,0){\bf \large {1A}}}
\put(7.5,8.5){\makebox(0,0){\bf \large {1B}}}
\put(2.0,7.0){\oval(5.0,1.0)\makebox(0,0){\bf \large {APC+Antigen}}}
\put(2.0,4.0){\oval(5.0,1.0)\makebox(0,0){\bf \large {Idiotopic cell}}}
\put(2.0,1.0){\oval(5.0,1.0)\makebox(0,0){\bf \large {Anti Id cell}}}

\put(1.5,3.5){\vector(0,-1){2}}
\multiput(2.5,1.5)(0,0.2){10}{\circle*{0.05}}
\put(2.5,3.5){\vector(0,1){0}}

\put(2.0,6.5){\vector(0,-1){2}}

\put(11.0,10.0){\oval(5.0,1.0)\makebox(0,0){\bf \large {APC+Antigen}}}
\put(11.0,7.0){\oval(5.0,1.0)\makebox(0,0){\bf \large {Naive Id cell}}}
\put(7.5,4.0){\oval(4.0,1.0)\makebox(0,0){\bf \large {Th1 Id cell}}}
\put(14.5,4.0){\oval(4.0,1.0)\makebox(0,0){\bf \large {Th2 Id cell}}}
\put(7.5,1.0){\oval(4.0,1.0)\makebox(0,0){\bf \large {Th1 cytokine}}}
\put(14.5,1.0){\oval(4.0,1.0)\makebox(0,0){\bf \large {Th2 cytokine}}}

\put(11,9.5){\vector(0,-1){2}}

\put(6.5,3.5){\vector(0,-1){2}}
\put(15.5,3.5){\vector(0,-1){2}}
\put(14.5,1.5){\vector(0,1){2}}
\put(7.5,1.5){\vector(0,1){2}}

\put(13.5,7){\line(1,0){2}}
\put(8.5,7){\line(-1,0){2}}
\put(15.5,7){\vector(0,-1){2.5}}
\put(6.5,7){\vector(0,-1){2.5}}
\put(14.5,4.5){\vector(0,1){2.5}}
\put(7.5,4.5){\vector(0,1){2.5}}

\multiput(8,1.5)(0.22,0.2){28}{\circle*{0.05}}
\multiput(14,1.5)(-0.22,0.2){28}{\circle*{0.05}}
\put(14.0,7.0){\vector(1,1){0}}
\put(8.0,7.0){\vector(-1,1){0}}
\end{picture}
\caption{The basic features of the model. 1A-Regulation of the id cell population by anti-id cells. 1B - Competition between TH1 T cells and TH2 T cells. Full lines designate activation. Broken lines designate suppression.}
\end{figure}

\begin{figure}
\center
\noindent
\setlength{\unitlength}{0.9cm}
\begin{picture}(20.0,15.0)

\put(10.0,14.0){\oval(5.0,1.0)\makebox(0,0){\bf \large {APC+Antigen}}}
\put(10.0,11.0){\oval(5.0,1.0)\makebox(0,0){\bf \large {Naive Id}}}
\put(6.0,8.0){\oval(3.0,1.0)\makebox(0,0){\bf \large {TH1 Id}}}
\put(14.0,8.0){\oval(3.0,1.0)\makebox(0,0){\bf \large {TH2 Id}}}
\put(6.0,5.0){\oval(3.0,1.0)\makebox(0,0){\bf \large {TH1 cytokine}}}
\put(14.0,5.0){\oval(3.0,1.0)\makebox(0,0){\bf \large {TH2 cytokine}}}
\put(2.0,8.0){\oval(3.0,1.0)\makebox(0,0){\bf \large {C cytokine}}}
\put(10.0,3.0){\oval(4.0,1.0)\makebox(0,0){\bf \large {Anti Id}}}
\put(2.0,5.0){\oval(3.0,1.0)\makebox(0,0){\bf \large {M$\Phi$}}}
\put(10.0,13.5){\vector(0,-1){2}}

\multiput(10.5,3.5)(-0.0,0.2){36}{\circle*{0.05}}
\put(10.5,10.5){\vector(0,1){0}}
\put(9.5,10.5){\vector(0,-1){7}}

\multiput(7.0,7.5)(0.06,-0.2){21}{\circle*{0.05}}
\put(7.0,7.5){\vector(0,1){0}}
\put(7.5,8.0){\vector(1,-4){1.1}}

\multiput(13.0,7.5)(-0.06,-0.2){21}{\circle*{0.05}}
\put(13.0,7.5){\vector(0,1){0}}
\put(12.5,8.0){\vector(-1,-4){1.1}}

\put(5.5,7.5){\vector(0,-1){2}}
\put(14.5,7.5){\vector(0,-1){2}}
\put(13.5,5.5){\vector(0,1){2}}
\put(6.5,5.5){\vector(0,1){2}}

\put(12.5,11){\line(1,0){2}}
\put(7.5,11){\line(-1,0){2}}
\put(14.5,11){\vector(0,-1){2.5}}
\put(5.5,11){\vector(0,-1){2.5}}
\put(13.5,8.5){\vector(0,1){2.5}}
\put(6.5,8.5){\vector(0,1){2.5}}

\multiput(14,4.5)(-0.26,-0.2){8}{\circle*{0.05}}
\put(12.0,3.0){\vector(-1,-1){0}}

\put(4.5,5.0){\vector(-1,0){1}}
\put(1.5,5.5){\vector(0,1){2}}
\put(2.5,7.5){\vector(0,-1){2}}
\put(3.5,8.0){\vector(1,0){1}}

\bezier{12}(8.0,3.0)(8.0,1.)(9.5,2.5)
\put(9.5,2.5){\vector(1,1){0}}
\end{picture}
\caption{Flow chart of the model. The model is based on a TH1/TH2 balance. If the TH1 concentration is too high the target organ is attacked. If there is a high TH2 concentration, the target organ remains healthy. Dotted lines indicate inhibition; full lines indicate stimulation. There is also an inhibition of TH2 on the production of TH1 and vice versa, but these lines were not drawn in order not to overload the graph. }
\end{figure}
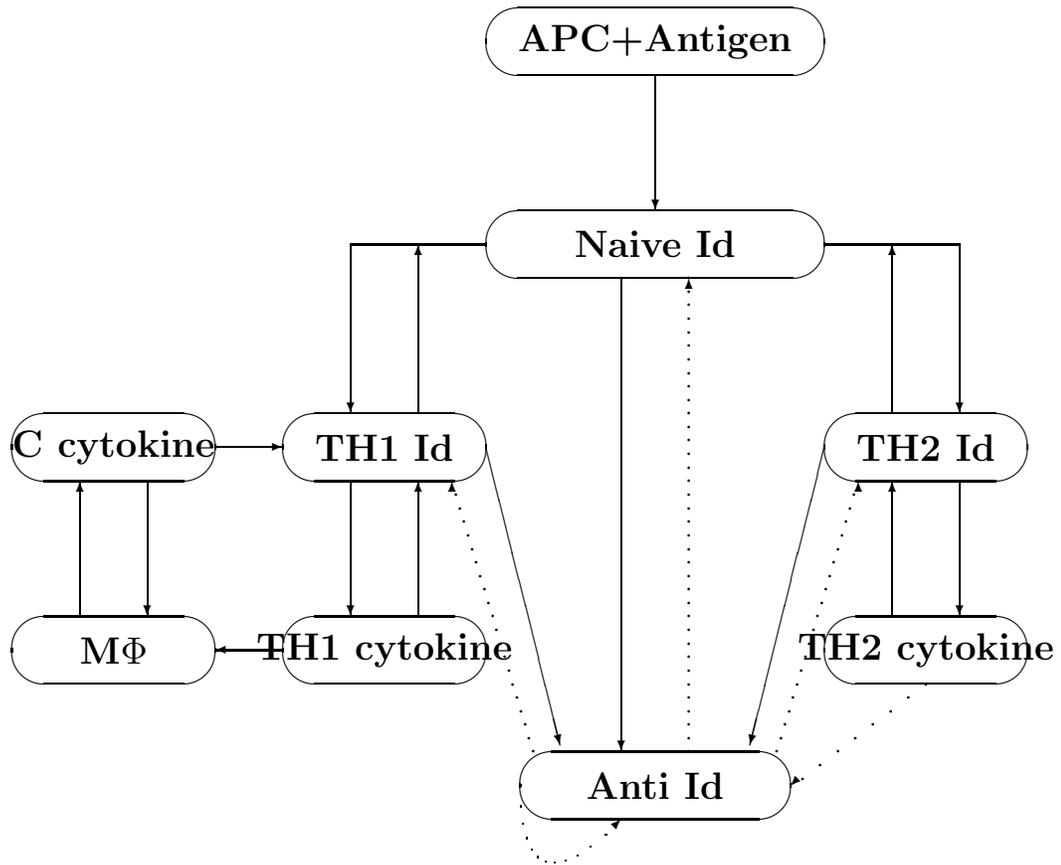

\subsection{Model mechanism}

We have two basic steady states (Figure 3), and all the cells and cytokines administrated only move the model from one steady state to the other, and back. 
The basic situation is a TH2 type steady state, which is assumed to be healthy, whereas the TH1 steady state leads to autoimmune disease. In contrast to the general view, the TH2 steady state represents a global situation of the system, and not the concentration of a single cell or cytokine. Note that the absence of disease does not mean that the system is at rest; the TH2 steady state is associated with relatively low concentrations of reacting cells. Nevertheless the id and ant-id cells are still reacting at the same rate, but at a lower concentration. In other words the healthy state merely looks like it is at rest; on the contrary health like disease is a state of interactions \cite{R30}.
\begin{figure}
\center
\noindent
\includegraphics[clip,width=11 cm]{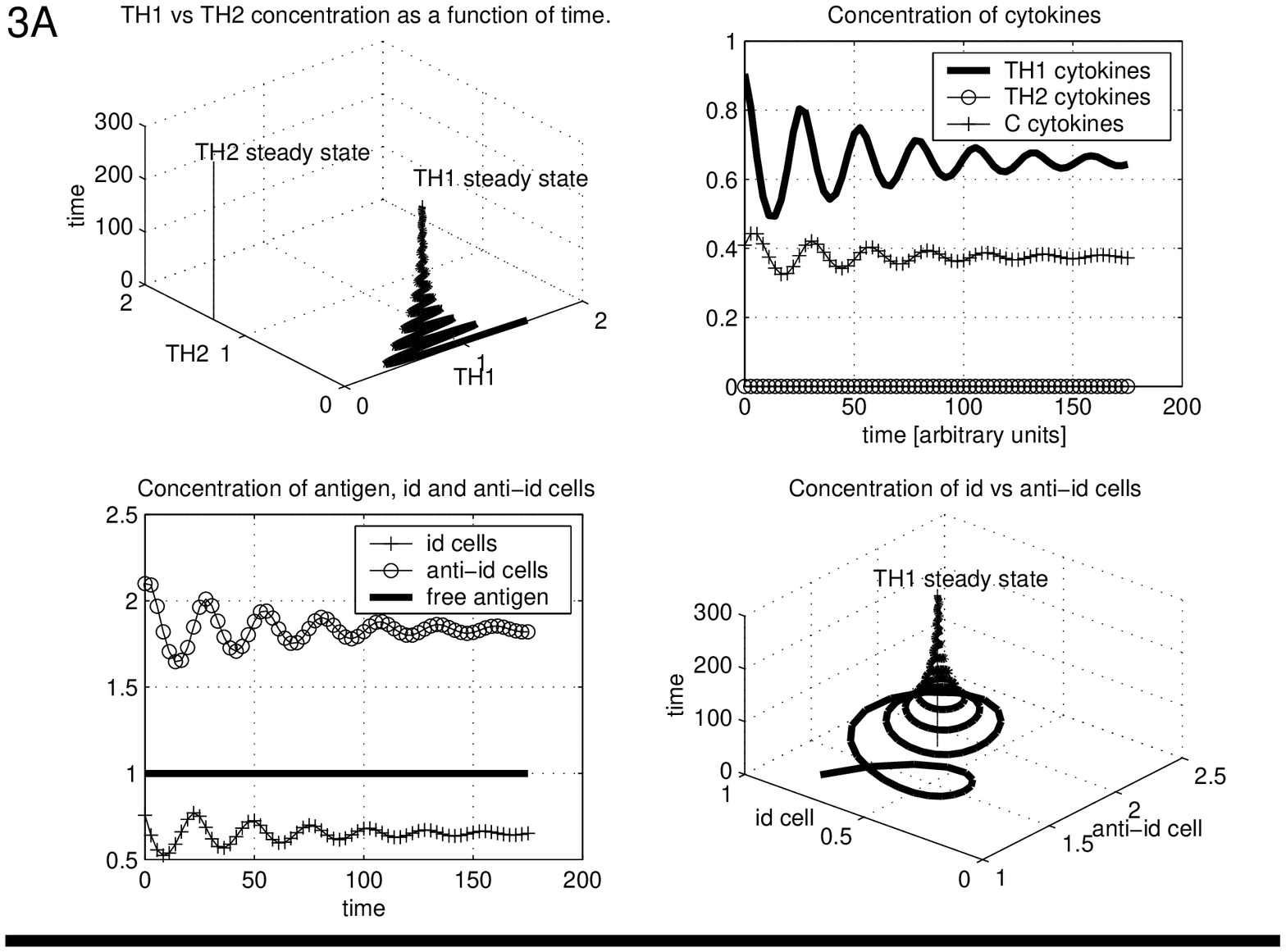}
\includegraphics[clip,width=11 cm]{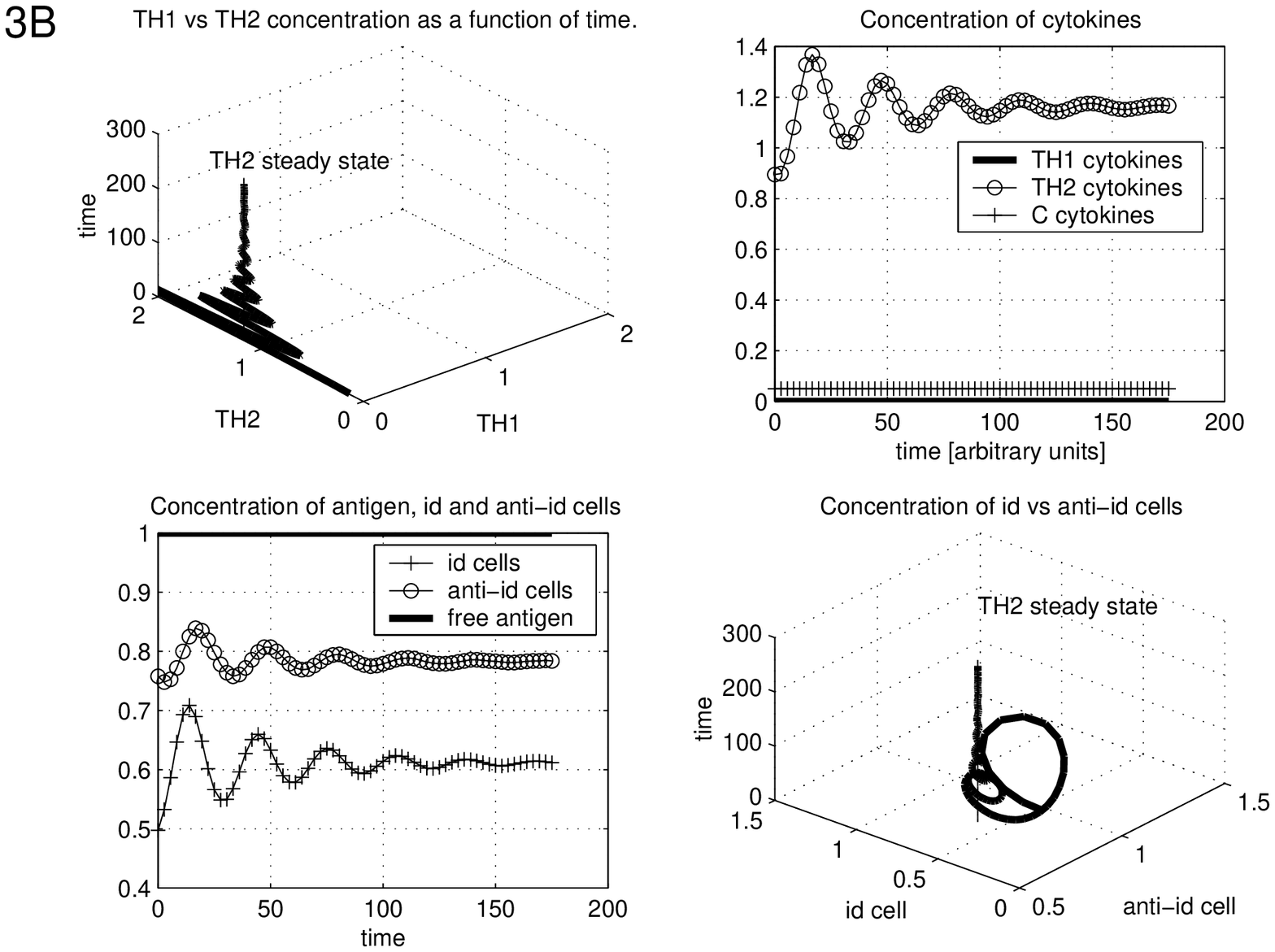}
\caption{The basic steady states. The first steady state is TH1 dominated (3A) and the second is TH2 dominated (3B). The model is such that, starting from a naive state, the system equilibrates into a TH2 steady state. 3A and 3B show concentrations in arbitrary units. Both 3A and 3B are composed of 4 graphs. The evolution in time of the TH1 and TH2 T cells concentration is depicted in the upper left drawing. The upper right drawing shows the cytokines concentration as a function of time. The drawing in the second row to the left shows the average concentration of id cells (TH1, TH2 and naive cells), the concentration of anti-id cells and the concentration of free antigen. The drawing at the right of the second row show the evolution in time of the id and the anti-id cell concentration. Note that the TH2 steady state (3B) manifests a lower concentration of both id and anti-id cells.}
\end{figure}

When we administer cells or antigen and adjuvant, we push the system toward the TH1  state (Figure 4). 
Following administration of T cells, which are primarily TH1 type cells, the number of TH1 cells is high enough to start the production of TH1 cytokines that inhibit the TH2 cytokines. 
The system then falls into the TH1 steady state, which is the diseased state. Another way to induce disease is to immunize the animal with antigen and adjuvant. We propose that the adjuvant stimulates macrophages to produce C type cytokines. The antigen together with the cytokines induces CD4 T cell proliferation and differentiation into TH1 cells, bringing the system again to a TH1 type steady state.
\begin{figure}
\center
\noindent
\includegraphics[clip,width = 11 cm]{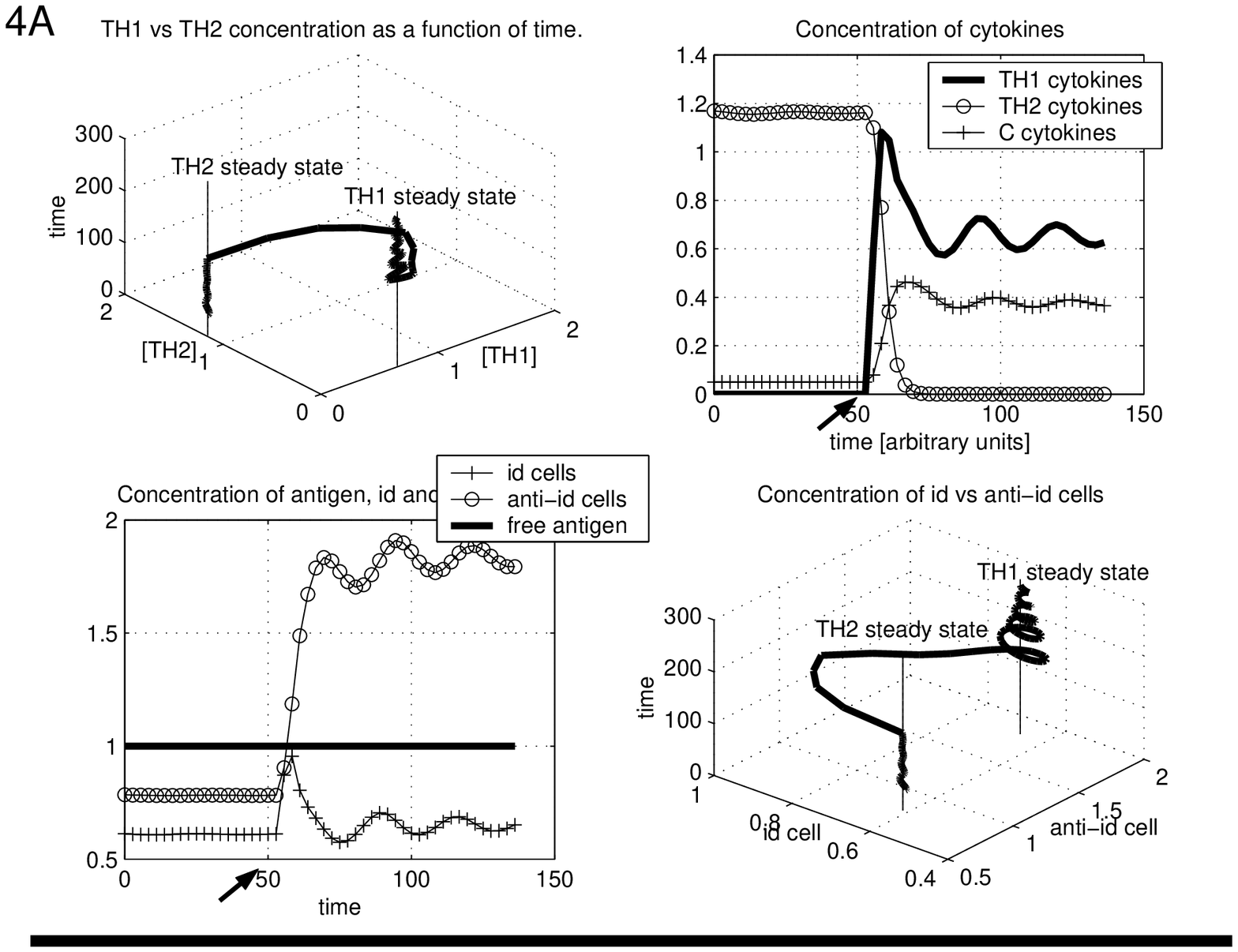}
\includegraphics[clip,width = 11 cm]{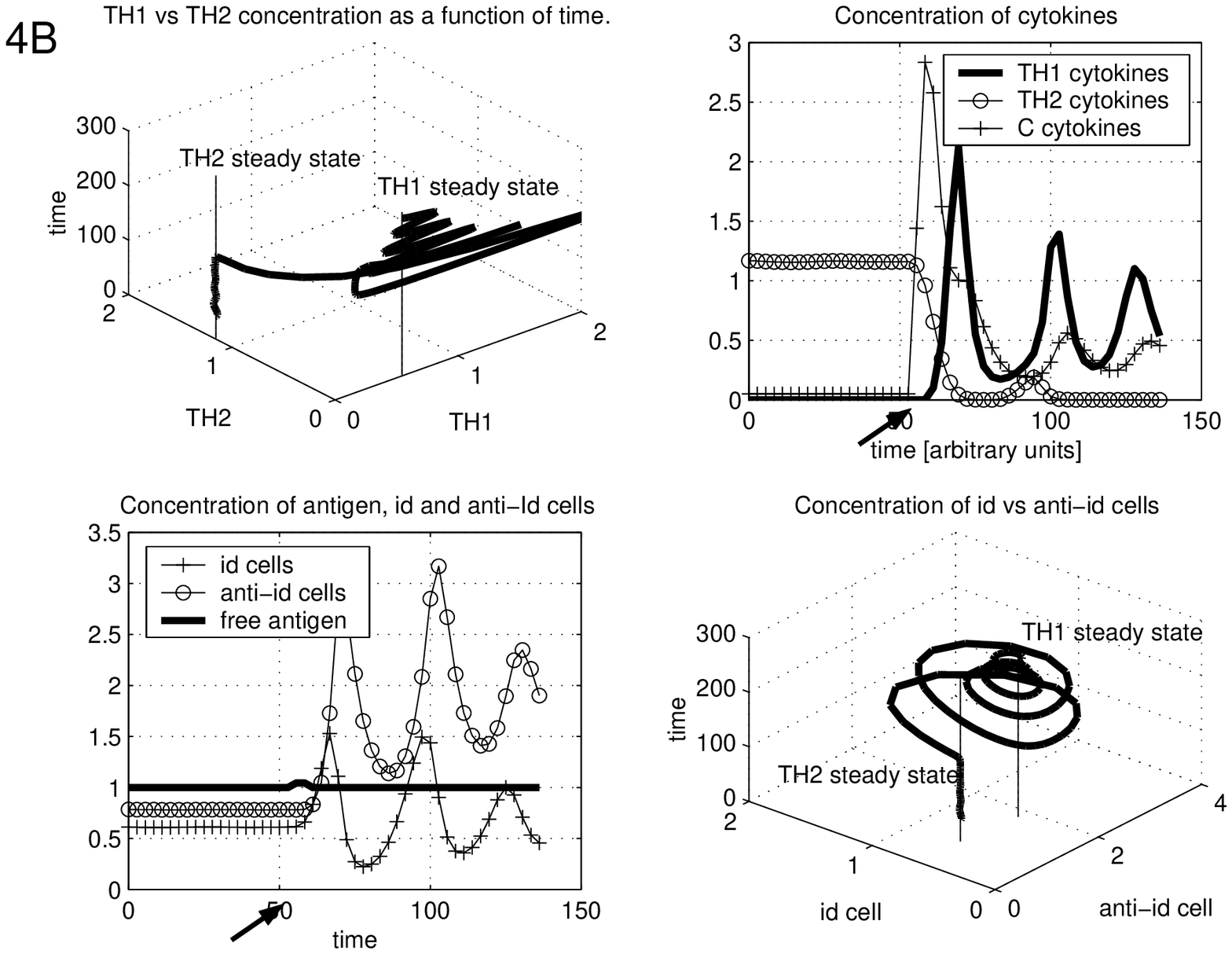}
\caption{The Transition from TH2 (healthy) steady state to a TH1 (disease) steady state. 4A shows the effect of adminstrating (arrow) TH1 cells. 4B shows the effect of immunization with antigen and adjuvant (arrow). Note that the TH1 steady state is associated with enhancement of the anti-id cells.}
\end{figure}

One way to heal the disease is to administer free antigen (Figure 5) in an amount that can raise the total number of naive CD4 cells that become TH1 CD4 cells. The resulting high concentration of CD4 cells will activate an increase in the number of anti-id cells. 
At this point, we assume that the lifespan of the anti-id cells is longer than that of the CD4 cells; the number of anti-id cells will remain high even after the number of CD4 id cells decreases (due to the fall in antigenic stimulation). Once the total number of CD4 cells is below a certain threshold, the system  returns by itself to the healthy TH2 situation.The basic situation is the TH2 situation.

\begin{figure}
\center
\noindent
\includegraphics[clip,width=11 cm]{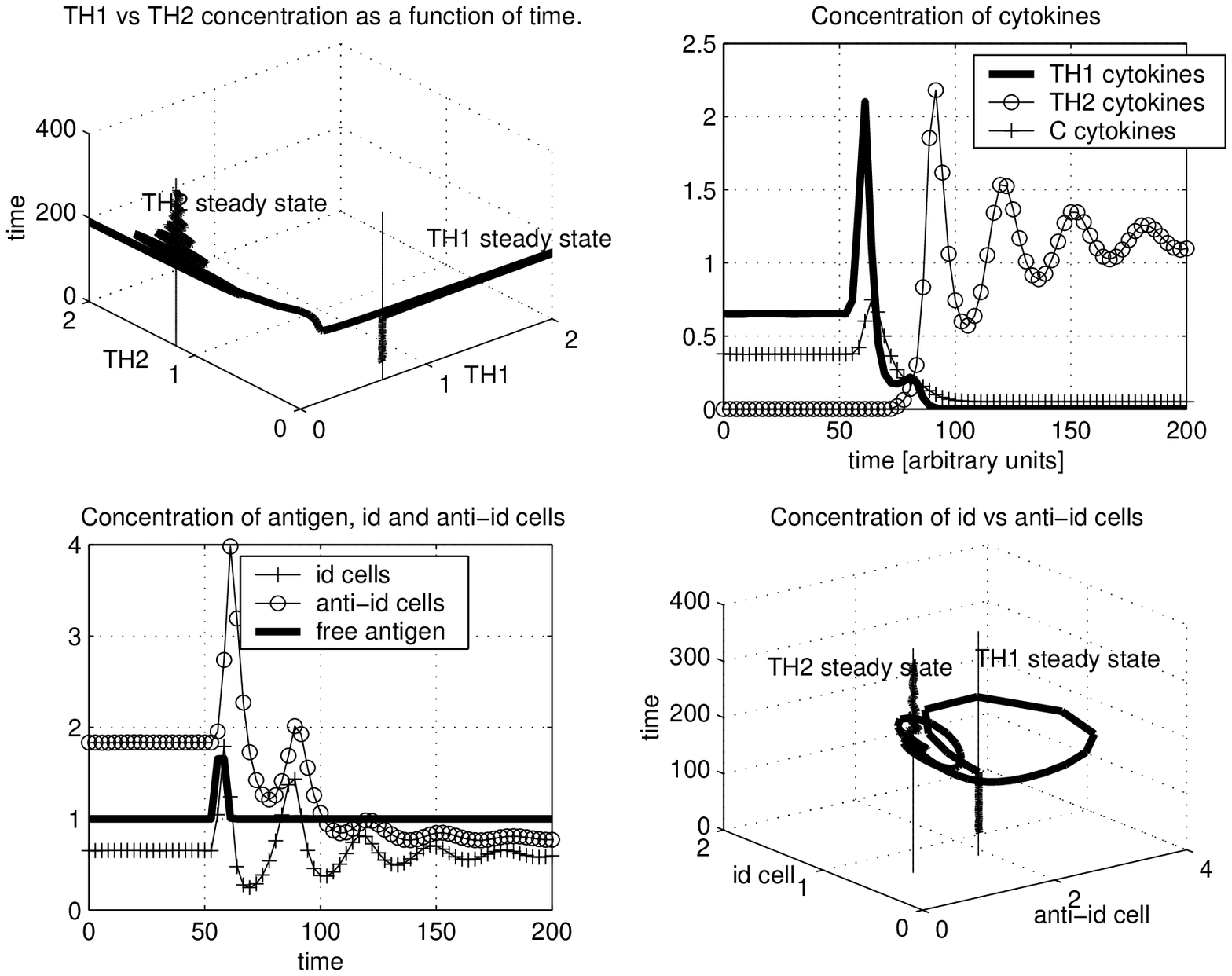}
\caption{ Healing with free antigen. Administrating free antigen raises the TH1 cell concentration, and then stimulates anti-id cells. Folowing it the anti-id cells concentration rises. The anti-id in turn supress the system and it returns back to its natural TH2 steady state. The notation used in this figure is identical to the notation used in figure 3.}
\end{figure}

The return to a TH2 state, following antigen administration, can be paradoxically inhibited by the administration of IL-4 with the antigen (Figure 6). The reason is as follows:  When we administer IL-4 (which is represented by type B cytokines), we inhibit the production of anti-id cells. In the absence of a sufficient number of anti-id cells, the number of CD4 effectors cells  remains above the threshold and the system stays in a TH1 state. 
If, however, the quantity of IL-4 is high enough to overcome TH1 dominance and starts activating TH2 T cells directly, then the TH1 cell concentration will decrease due to the production of TH2 cytokines by the newly formed TH2 cells, despite the low number of anti-id cells. 
\begin{figure}
\center
\noindent
\includegraphics[clip,width=11 cm]{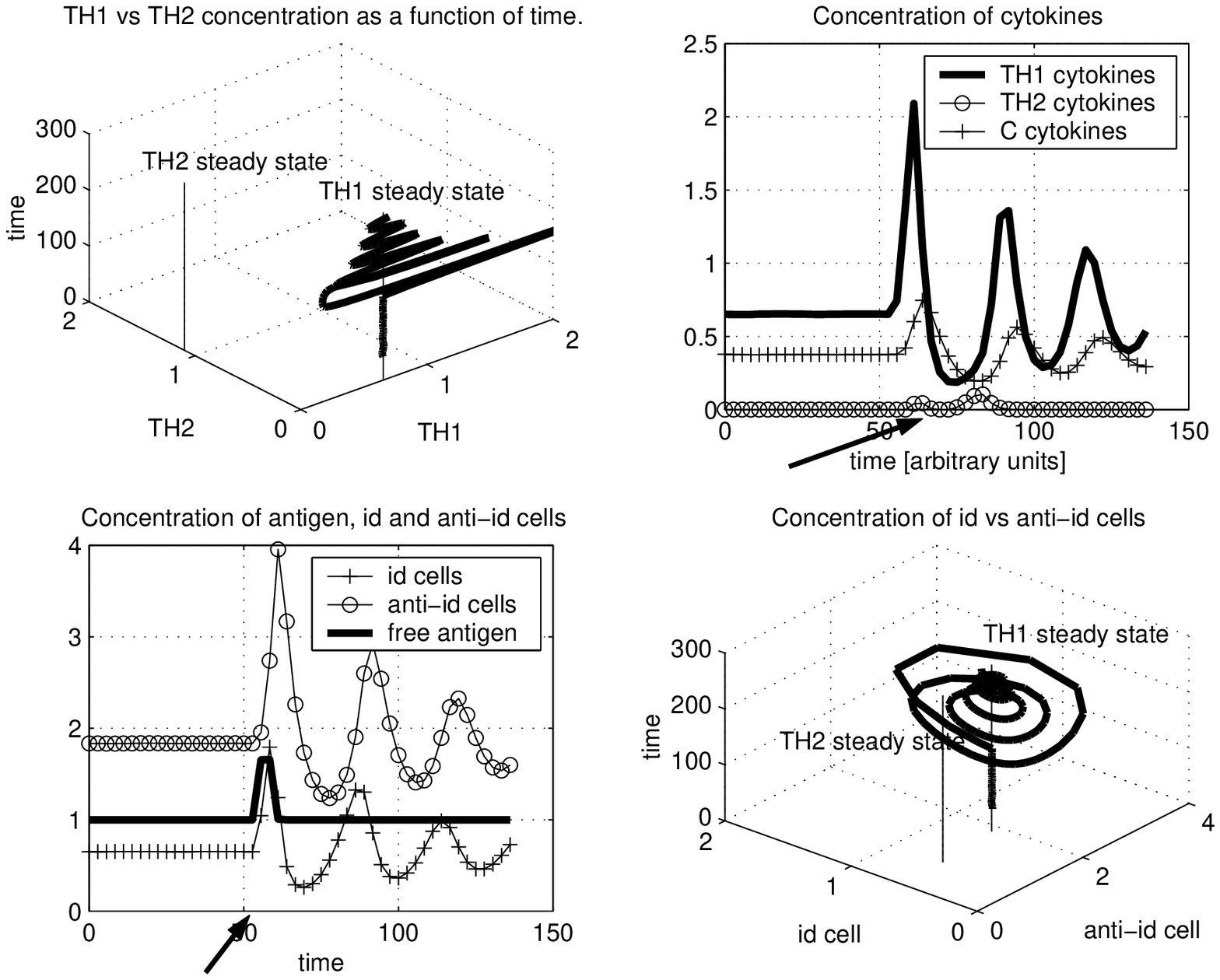}
\caption{Inhibition of the healing effect of free antigen by administration of IL4. IL4 has a double role in this case. At the first stage it reduces the number of TH1 id cells, and at the second stage it reduces the number of anti-id cells, so that their effect is weakened. The notation used in this figure is identical to the notation used in figure 3.}
\end{figure}

The administration of anti-IL-4 antibody with the antigen can paradoxically act in a manner similar to the administration of IL-4 with the antigen. Note that in both cases the inhibition of disease is weaker than the effect of administrating free antigen.
The administration of anti-IL-4 antibody has no effect on the system at first, since the concentration of IL-4 is low. 
The effect of the antigen is then as described at the beginning of the section;
the antigen induces a sharp rise in the concentration of id CD4 cells, and following the id CD4 cells, the anti-id cell concentration raises. These anti-id cells produce a long term decrease in the concentration of the CD4 id cells. 
At a later stage, when the concentration of CD4 id cells has decreased to a low level, the effect of the anti-IL-4 antibody becomes important. 
The level of id CD4 TH1 type cells at this stage is low, and the system should now start to produce TH2 type cells and cytokines. However the anti-IL-4 antibody now inhibits the production of IL-4, and the system has a lower probability to return to a TH2 state.

Just as the administration of TH2 cytokines can inhibit the transition back to a TH2 state, the administration of TH1 cytokines can rise the probability of such a transition.
The continuous administration of TH1 cytokines will of course keep the cytokines in the TH1 profile and induce or enhance the disease. 
Administration of the TH1 cytokines at the onset of the disease (Figure 7) will lead to an increase in the production rate of TH1 id CD4 cells; 
the high concentration of TH1 cells will raise the concentration of the anti-id cells, which will later feed back to decrease the concentration of CD4 id cells. 
If the CD4 id cell concentration decreases to a low level, the system will return to its natural steady state, which is a TH2 steady state, as described above.

The administration of some cytokines has been shown to have no observable effect on the system. This does not mean that these cytokines have no important role in the dynamics of this system. In fact the lack of an observable effect can be due to a high endogenous concentration of the cytokines; adding additional amounts cannot enhance the effect.

\begin{figure}
\center
\noindent
\includegraphics[clip,width=11 cm]{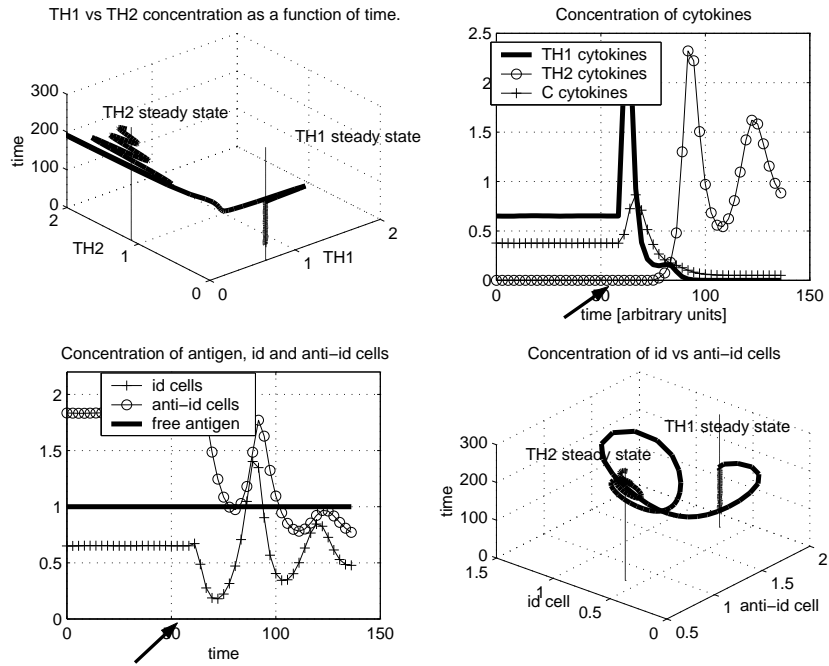}
\caption{ Healing by administrating TH1 cytokines. The TH1 Cytokines raise the TH1 cell concentration, inducing a rise in the concentration of anti-id cells. The anti-id cells lower the concentration of the id cells, returning the system back to its normal (healthy) TH2 steady state. The notation used in this figure is identical to the notation used in figure 3.}
\end{figure}

Administration of a cytokine may have no effect if the cytokine concentration is close to, or above its threshold for causing an effect. For example, administration of IL-12 at the onset of disease may have no effect if the IL-12 effect is already at its maximum. Paradoxically in this case, the administration of a low dose of anti-IL-12 antibody will also have no effect, since the precise concentration of the molecule is not important. A large change in the concentration of the molecule leads to a small change in its effects. This can be understood if the effect is bell shaped as a function of the concentration. At low concentrations, a small change in the concentration leads to a large change in the function value. However at the top of the bell, a large change in the concentration leads to very small change in the function value. 

\subsection{Redundancy and gene knock-out mice}
Gene knock-out mice lacking some of the TH1 cytokines can still develop a TH1-dependent disease. However this does not contradict the model we presented, which describes the evolution of the disease as the transition to a TH1 type cytokine profile. The cytokines in this model are described as a group with the same approximate action. Even if such a group lacks one of its components, the effect of the group as a whole is not altered.
The TH1/TH2 phenotype does not depend on a single cytokine, since the concentration of each one of the TH1 and TH2 cytokines can change by a large number of mechanisms. The immune system cannot afford  to change its general phenotype each time an external factor modifies the concentration of a single cytokine. The immune system achieves its robustness to external changes in the cytokine profile through redundancy \cite{R30}. A large number of cytokines can each have common effects. Thus a certain result can occur even if not all the cytokines leading to the effect are expressed. 

Indeed, the fact that the body uses many redundant cytokines can have a larger meaning. Consider that the effects mediated by cytokines have at least two functions. The first function is to maintain a healthy steady state, and the second function relates to the fine details of the steady state. 
In normal behavior, the changes are felt at the level of the details, but in extreme situations, like the response to an external pathogen, the body will try to change its total steady state and produce a fast and strong response.

\section{Discussion}
\subsection{Regulation of autoimmunity}
This paper presents a simplified model of the kinetics of the concentrations of the main cell and cytokine types in the evolution of a TH1 type autoimmune disease. It shows that even simple dynamics can explain many of the complex features in the  evolution of the response. 
We have grouped the major cells and cytokines into single types, although there are obviously differences in the roles taken by each one of the cells we have grouped together.
 This model ignores the effects of localization and of tissue barriers (for example, the blood/brain barrier in EAE) in inhibiting the passage of some of cells and allowing others.
We did not take into account the effects of the tissue cells (for example astrocytes in EAE, or  islet cells in IDDM), of the cytokines they secrete, and of B cells and antibodies. Some major mechanisms were ignored including anergy and the release of free antigen due to cell destruction.

All these simplifications render this a general feature model rather than a precise quantitative description of the cell concentrations. 
Our model is built on a double system of positive and negative feedback between id and anti-id cells, and between TH1 and TH2 states. 
Such a structure could be the basic feature of the selection between TH1 and TH2 for a group of important antigens in the body. 
This model seems to describe the situation both in EAE and in IDDM, and we have no reason to believe that it is not general. Negative feedback keeps the general concentrations of all kinds of cells within reasonable limits, and positive feedback enables the system to mount a fast enough response to an antigen when needed. 
With this structure, it is possible to maintain a population of T cells responding to self-antigens of the body, keeping them in a TH2 state where they may be restrained.
\subsection{Models}
This kind of model shows the complex dynamics attained even by a small system containing only 5 cell types and 3 cytokine types. It shows that even small systems need a kinetic description, and the system cannot be understood by an overall estimation of the effect of one cell type on the others. 
Such models must be robust to changes in parameters, since it is obvious that the interactions between the cells in the body are not precisely fixed and that the cells are constantly subjected to chemical noise impinging on them. Nevertheless, the body must keep its general steady state, although the precise quantity of any cell is not very important.
\subsection{Feedback loops}
The present type of modeling shifts our attention from the numbers of each cell type to the global state of the immune system, which most of the time should be in some quasi-equilibrium between many opposing influences. 
The way to analyze such a global state is by positive and negative feedback loops and the interactions between them. 
According to this formalism, our model can be described as a large negative feedback loop, built of 4 internal loops :
\begin{itemize}
\item The id-anti-id cell negative feedback loop, in which are nested three positive and one negative feedback loop. These include:
\item The TH1 and TH2 cells each of which feedsback positively on itself.
\item The TH1 and macrophage positive feedback.
\item The TH2 and anti-id cell negative feedback. 
\end{itemize}
In conclusion, the description of disease evolution can be based on a global kinetic view, and the subsystems are not cell types but feedback loops. (figure 8)

\begin{figure}
\center
\noindent
\setlength{\unitlength}{0.9cm}
\begin{picture}(8.0,21.0)
\thicklines
\put(4.0,20.0){\makebox(0,0){\bf \large {8A}}}
\put(2.0,15.0){\makebox(0,0){\bf \large {8B}}}
\put(6.0,15.0){\makebox(0,0){\bf \large {8C}}}
\put(4.0,10.0){\makebox(0,0){\bf \large {8D}}}
\put(4.0,5.0){\makebox(0,0){\bf \large {8E}}}

\put(4.0,19.0){\oval(4.0,1.0)\makebox(0,0){\bf \large {id cell}}}
\put(4.0,16.0){\oval(4.0,1.0)\makebox(0,0){\bf \large {Anti-id cell}}}

\put(3.0,18.5){\vector(0,-1){2}}
\multiput(5.0,16.5)(0,0.2){10}{\circle*{0.05}}
\put(5.0,18.5){\vector(0,1){0}}

\put(2.0,14.0){\oval(3.8,1.0)\makebox(0,0){\bf \large {TH1 id cell}}}
\put(6.0,14.0){\oval(3.8,1.0)\makebox(0,0){\bf \large {TH2 id cell}}}
\put(2.0,11.0){\oval(3.8,1.0)\makebox(0,0){\bf \large {TH1 cytokine}}}
\put(6.0,11.0){\oval(3.8,1.0)\makebox(0,0){\bf \large {TH2 cytokine}}}

\put(1.0,13.5){\vector(0,-1){2}}
\put(5.0,13.5){\vector(0,-1){2}}
\put(3.0,11.5){\vector(0,1){2}}
\put(7.0,11.5){\vector(0,1){2}}

\put(2.0,9.0){\oval(2.0,1.0)\makebox(0,0){\bf \large {C cyt}}}
\put(6.0,9.0){\oval(2.0,1.0)\makebox(0,0){\bf \large {TH1 id}}}
\put(2.0,6.0){\oval(2.0,1.0)\makebox(0,0){\bf \large {M$\Phi$}}}
\put(6.0,6.0){\oval(2.0,1.0)\makebox(0,0){\bf \large {TH1 cyt}}}
\put(3.0,9.0){\vector(1,0){2.0}}
\put(5.0,6.0){\vector(-1,0){2.0}}
\put(2.0,8.5){\vector(0,-1){2}}
\put(6.0,6.5){\vector(0,1){2}}

\put(6.0,4.0){\oval(2.0,1.0)\makebox(0,0){\bf \large {TH2 id}}}
\put(2.0,2.5){\oval(2.0,1.0)\makebox(0,0){\bf \large {anti-id}}}
\put(6.0,1.0){\oval(2.0,1.0)\makebox(0,0){\bf \large {TH2 cyt}}}
\put(5.1,3.5){\vector(-2,-1){2}}
\put(6.5,3.5){\vector(0,-1){2}}
\put(5.5,1.5){\vector(0,1){2}}
\multiput(2.0,3.0)(0.2,0.07){15}{\circle*{0.05}}
\multiput(2.0,2.0)(0.2,-0.07){15}{\circle*{0.05}}
\put(5.0,4.1){\vector(1,1){0}}
\put(2.05,2.0){\vector(-1,1){0}}

\end{picture}
\caption{ A schematic description of the five feedback loops. 8A Shows the negative feedback loop between id and the anti-id cells. 8B and 8C show the positive feedback of the TH1 and TH2 cell and cytokines on themselves. 8D shows the positive feedback of the macrophages and the TH1 id cells and macrophages. 8E shows the negative feedback between the TH2 cytokines and the anti-id cells.}
\end{figure}
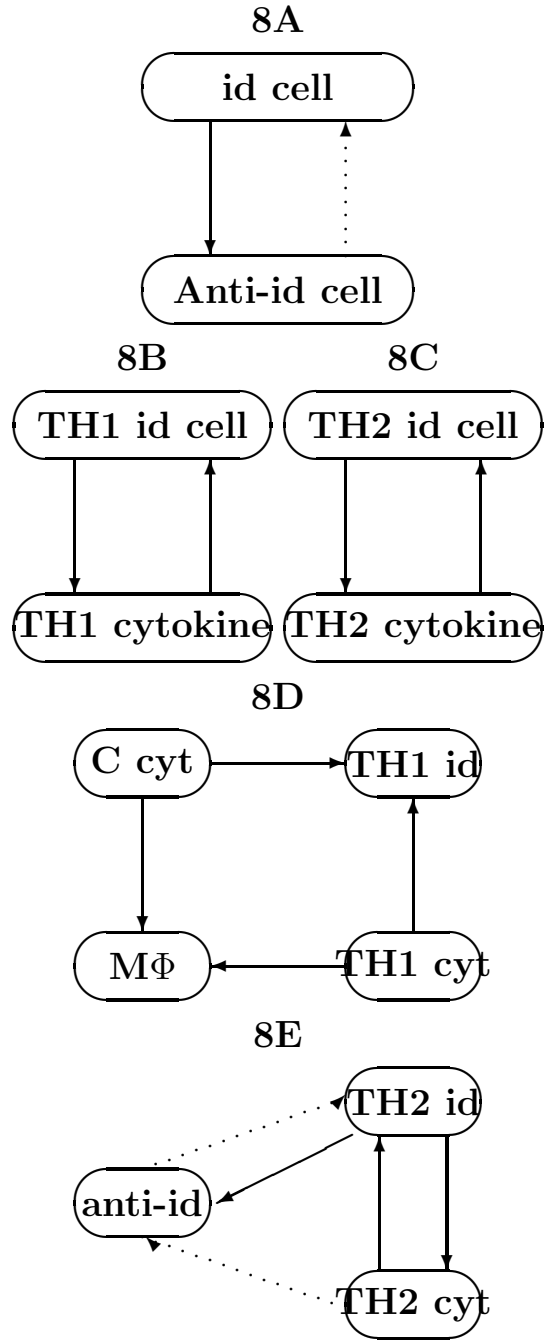

\section{Acknowledgment}
We wish to thank Prof G. Steinman for the useful data on EAE he graciously supplied.

Irun R.Cohen is the incumbent of the Mauerberg Chair in Immunology, Director of the Robrt Kock-Minerva Center for Research in Autoimmune Disease, and Director of the Center for the Study of emerging Diseases.

\end{document}